\documentclass[draft,numreferences,thms]{kluwer}
\usepackage{amsopn}
\usepackage{amstext}
\usepackage{amscd}
\newcommand{\To}{\longrightarrow}
\newcommand{\n}{\mathcal{N}}

\newcommand{\cc}{{\mathcal C}}

\newcommand{\CC}{C_{abcd}}
\newcommand{\CCC}{C_{a'b'c'd'}}
\newcommand{\prl}{\mathbb{P}^{1}}
\newcommand{\abs}[1]{\left|#1\right |}

\newcommand{\tq}{\,|\,}
\newcommand{\md}[1]{\ \mbox{\rm(mod }{#1})}
\newcommand{\h}[1]{H^1(G_k,{#1})}
\newcommand{\gen}[1]{\langle\, {#1} \,\rangle}
\newcommand{\io}{\gen{\iota}}
\newcommand{\op}{\operatorname}
\newcommand{\as}{\op{AS}}
\newcommand{\kas}{k/\as(k)}
\newcommand{\ke}{\op{Ker}(E)}

\newcommand{\ck}{\op{Cub}_k}
\newcommand{\au}{\op{Aut}(C)}
\newcommand{\tw}{\op{Tw}(C/k)}
\newcommand{\twc}{\op{Tw}^0(C/k)}
\newcommand{\g}{\Gamma}
\newcommand{\gw}{\Gamma_W}
\newcommand{\gu}{\Gamma_{u(x)}}
\newcommand{\gab}{\Gamma_{abcd}}
\newcommand{\ga}{\gamma}

\newcommand{\la}{\lambda}
\newcommand{\ff}[1]{\mathbb F_{q^{#1}}}
\newcommand{\f}[1]{\mathbb{F}_{#1}}
\newcommand{\fs}{\mathbb{F}_{16}}
\newcommand{\pg}{\op{PGL}_2(k)}
\newcommand{\kb}{\overline{k}}
\newenvironment{pmat}{%
  \left(\begin{array}{cccc}}{\end{array}\right)}
\newenvironment{pmat1}{%
  \left(\begin{array}{c}}{\end{array}\right)}
\newenvironment{pmat3}{%
  \left(\begin{array}{ccc}}{\end{array}\right)}
  {\end{eqnarray}%
}
\newenvironment{myeqnarray*}
{\setlength{\arraycolsep}{1.5pt}\begin{eqnarray*}}%
{\end{eqnarray*}}%

\begin{document}

\begin{article}

\begin{opening}

  \title{Curves of genus two over fields of even characteristic}

  \author{Gabriel \surname{Cardona}
    \email{gabriel.cardona@uib.es}
    \thanks{supported by BFM2000-0794-C02-02 and HPRN-CT-2000-00114}
  }
  \institute{
    Dept. Ci\`encies Matem\`atiques i Inf.\\
    Universitat de les Illes Balears\\ Ed. Anselm Turmeda, Campus UIB\\
    Carretera Valldemosa, km. 7.5\\ E-07071 Palma de Mallorca, Spain
  }

  \author{Enric \surname{Nart}
    \email{nart@mat.uab.es}
  }
  \author{\ \\Jordi \surname{Pujol{\`a}s}
    \email{pujolas@mat.uab.es}
    \thanks{supported by BHA2000-0180}
  }
  \institute{
    Departament de Matem\`atiques\\
    Universitat Aut\`onoma de Barcelona\\ Edifici C\\ E-08193 Bellaterra,
    Barcelona, Spain
  }

  \runningtitle{Curves of genus two over fields of even
  characteristic}
  \runningauthor{G.~Cardona, E.~Nart and J.~Pujol{\`a}s}
  
  \begin{ao}\\
    Gabriel Cardona\\
    Dept. Ci\`encies Matem\`atiques i Inf.\\
    Universitat de les Illes Balears\\ Ed. Anselm Turmeda, Campus UIB\\
    Carretera Valldemosa, km. 7.5\\ E-07071 Palma de Mallorca, Spain\\
  \end{ao}

  \begin{abstract}
    In this paper we classify curves of genus two over a perfect field
    $k$ of characteristic two. We find rational models of curves with a
    given arithmetic structure for the ramification divisor and we give
    necessary and sufficient conditions for two models of the same type
    to be $k$-isomorphic. As a consequence, we obtain an explicit formula
    for the number of $k$-isomorphism classes
    of curves of genus two over a finite field. Moreover,
    we prove that the field of moduli of any curve coincides with its
    field of definition, by exhibiting rational models of curves with any
    prescribed value of their Igusa invariants. Finally,
    we use cohomological methods to find, for each
    rational model, an explicit description of its twists.
    In this way, we obtain a parameterization of
    all $k$-isomorphism classes of curves of genus two in terms of
    geometric and arithmetic invariants.
  \end{abstract}


  {Primary 11G20, 14G27; Secondary 14G50.} 
\end{opening}


\section*{Introduction}
Let $k$ be a perfect field of even characteristic. Igusa described in
\cite{Igusa} the moduli variety for curves of genus two
as certain 3-dimensional affine variety, whose $k$-points
are in bijection with $k^3$. Moreover, he gave explicit formulas to
compute the moduli point of a curve in terms of invariants. In
section 2 of this paper we prove that the field of moduli of any
curve coincides with its field of definition, so that the Igusa
invariants determine a bijection between
the set of
$\kb$-isomorphism classes of curves of genus two defined over $k$
and the affine space $k^3$.
This is known to be false for some fields of odd or zero
characteristic (cf. \cite{Mestre}). We prove this result by exhibiting
rational models of curves with any prescribed value of their
invariants.

Rational models of curves of genus two are studied more generally in
section 1, where we obtain an effective classification,
up to $k$-isomorphism,
of curves of
genus two defined over $k$. For each possible
arithmetic structure of the ramification divisor, the corresponding
curves admit quasi-affine models $y^2+y=u(x)$, where $u(x)\in k(x)$
is a rational function with a concrete divisor of poles. The
$k$-isomorphism classes of these curves are in bijection with the
orbits of these rational functions $u(x)$ under a  double
action by the Artin-Schreier group $\as(k(x))$ and the projective
linear group $\pg$. In paragraph 1.3 we carry out this classification
in a very explicit manner. As a by-product, we obtain an explicit
computation of the $k$-automorphism group of each
curve.

In section 3 we describe all twists of a given curve $C$, in terms of
the rational models of section 1. This is achieved by computing
$\h{\au}$ for any possible structure of $\au$. Together with the
result of section 2, this computation furnishes an explicit
parameterization of all $k$-isomorphism classes of curves of genus
two, each class being determined by a couple of invariants: one
geometric (a triple $(j_1,j_2,j_3)\in k^3$ of values of the Igusa
invariants) and the other arithmetic (a concrete twist of a chosen
curve having $(j_1,j_2,j_3)$ as geometric invariants).

Our initial motivation for this work was to classify curves of genus
two over finite fields of even characteristic, in view of their use
in Cryptography (cf. \cite{Kob}). In the finite field case our results
provide formulas for the number of curves of genus two with a fixed
structure of the ramification divisor (paragraph 1.4) and for the
number of curves with a fixed structure of the full automorphism
group (paragraph 3.2).

\section{Rational models of curves of genus two}
We fix once and for all a perfect field $k$ of even characteristic
and an algebraic closure $\kb$ of $k$. We denote by $G_k$ the Galois
group $\op{Gal}(\kb/k)$ and by $x\mapsto \sqrt x$ the automorphism of
$k$ inverse of the Frobenius automorphism, $x\mapsto x^2$. Also,
$\mathbb{F}_{2^m}$ and $\mu_n$ will denote respectively the unique
finite subfield of 
$\kb$ with $2^m$ elements and the subgroup of $\kb^*$ of the $n$-th
roots of unity. Finally, we denote $\mu_n(k)=\mu_n\cap k$.

By a {\it curve of genus two} we mean a smooth, projective,
geometrically irreducible curve of genus two.

\subsection{Hyperelliptic curves}
In this paragraph we review well-known results on hyperelliptic
curves. Let $C$ be a hyperelliptic curve defined over $k$; that is,
 $C$ is a smooth, projective and geometrically irreducible curve
defined over $k$, of genus $g\ge 2$, admitting a degree 2 morphism
$\pi\colon C \To \prl$, which is  also
defined over $k$.

The corresponding extension of function fields, $
k(C)/\pi^*(k(\prl))$, is a cyclic quadratic extension (if it were not
separable, then $k(C)$ would be a function field of genus zero). To
the non-trivial element of the Galois group it corresponds an
involution, $\iota\colon C\longrightarrow C$, which is called the
{\it hyperelliptic involution} of $C$. By definition, thus, two
points $P,Q\in C(\kb)$ have the same image under $\pi$ if and only if
$Q=P$ or $Q=P^{\iota}$.

We recall some basic properties of $\iota$ based on the crucial fact
that it is a canonical involution, independent of the morphism $\pi$.

\begin{thm}\label{restr}
  Let
  $\pi_1 ,\pi_2 \colon C \To \prl$ be two $k$-morphisms of degree $2$.
  Then, there exists a unique $k$-automorphism $\gamma$ of $\prl$ such
  that $\pi_2=\gamma\circ \pi_1$.
\end{thm}

\begin{cor}\label{hypinv} 
\begin{enumerate}
\item The hyperelliptic involution is canonical. The fixed points of
  $\iota$ are the Weierstrass points of $C$ and they are the
  ramification points of any morphism of degree 2 from $C$ to $\prl$.
\item Any $k$-automorphism $\varphi$ of $C$ fits into a commutative
  diagram:
  $$
  \begin{CD} C @>{\varphi}>> C \\ @V{\pi}VV @VV{\pi}V
    \\ \prl @>{\gamma}>> \prl\\
  \end{CD}
  $$
  for certain uniquely determined $k$-automorphism $\gamma \in
  \op{Aut}_k(\prl)$. The map $\ \varphi\mapsto \gamma\ $ is a
  homomorphism (depending on $\pi$) and we have an exact sequence of
  groups:
  $$
  1\To\io \To \op{Aut}_k(C) \stackrel{\pi}\To\op{Aut}_k(\prl).
  $$
  In particular, $\iota$ commutes with every automorphism of $C$.
\end{enumerate}
\end{cor}

We recall now some generalities on Artin-Schreier (i.e. cyclic
quadratic) extensions of a field $K$ of even characteristic. Let
$K^{\mbox{\tiny sep}}$ be a fixed separable closure of $K$ and
$G_K=\op{Gal}(K^{\mbox{\tiny sep}}/K)$. The Artin-Schreier
homomorphism, $x\mapsto\as(x)=x+x^2$, fits into an exact sequence:
$$
0\To\f2\To K^{\mbox{\tiny sep}}\stackrel{\as}\To K^{\mbox{\tiny
sep}}\To 0,
$$
leading to an isomorphism $H^1(G_K,\f2)\simeq K/\as(K)$.
Thus, the group $K/\as(K)$ classifies all cyclic quadratic
(or trivial) extensions of $K$.
Any $u\in K-\as(K)$ determines a separable quadratic extension of $K$
by adjoining to $K$ the roots of the separable irreducible polynomial
$Y^2+Y+u$. An element $u'\in K$ furnishes the same extension than $u$
if and only if $u+u'\in \as(K)$. Moreover, any cyclic quadratic
extension of $K$ can be obtained this way.
%
If $K$ is a finite field,
$\as(K)$ coincides with the subgroup of elements of
absolute trace zero and  $K/\as(K)$ has only two elements.

By fixing a point $\infty \in \prl (k)$, we fix
isomorphisms $k(\prl)\simeq
k(x)$ and $\op{Aut}_k(\prl)\simeq \pg$. The function field $k(C)$ is
then identified (via $\pi^*$) to an Artin-Schreier extension of
$k(x)$ and it admits an Artin-Schreier generator $y\in k(C)$
satisfying $\ y^2 + y = u(x)$, for certain $u(x) \in k(x)-\as(k(x))$.

Conversely, for any $u(x)\in k(x)-\as(k(x))$, the equation
$y^2+y=u(x)$ determines an Artin-Schreier extension, $L=k(x,y)$, of
$k(x)$. Such an equation determines a plane non-singular quasi-affine
curve $C^{\mbox{\tiny af}}=C^{\mbox{\tiny af}}_{u(x)}$, defined over
$k$.
We shall denote by $C=C_{u(x)}$ the projective, smooth curve
obtained as the normalization of the projective closure of
$C^{\mbox{\tiny af}}$. The projection on the first coordinate,
$(x,y)\mapsto x$, lifts to a morphism of degree 2, $\pi\colon
C\To\prl$, implicitly associated to the equation. Clearly, the
function field of $C$ is $k$-isomorphic to $L$ and the homomorphism
$\pi^*\colon k(\prl)\To k(C)$ translates into $k(x)\hookrightarrow
k(x,y)$ under natural identifications.

Since $C^{\mbox{\tiny af}}$ is non-singular, we can identify its
points with an open subset of $C(\kb)$; this allows us to attach
affine coordinates to most of the points of $C$. Hence, to deal with
these curves in practice, it is not necessary to find equations for
$C$ itself, as soon as one has control on the finite set of {\it
points at infinity}: $Z_{\infty}:=C(\kb)-C^{\mbox{\tiny af}}(\kb)$.
Note that the ramification points of $\pi$ lie above poles of $u(x)$;
hence, they are always points at infinity.

The genus $g$ of $C_{u(x)}$ is easily determined in terms of the
poles of odd order of the elements in $u(x)+\as(k(x))$, and $C_{u(x)}$
is always geometrically irreducible, hence a hyperelliptic curve,
when $g\ge 2$ (cf. Proposition \ref{ramdiv} below). In affine
coordinates, the hyperelliptic involution is then expressed as:
$(x,y)^{\iota}=(x,y+1)$. After Theorem \ref{restr}, it is easy to
determine in terms of $u(x)$ the $k$-isomorphism class of 
$C_{u(x)}$:

\begin{prop}\label{xy}
  Two hyperelliptic curves $C_{u(x)}$, $C_{u'(x)}$ are $k$-iso\-mor\-phic
  if and only if
  $$
  u'(x)\equiv u(\ga(x)) \md{\as(k(x))},
  $$
  for some $\ga\in\pg$.
\end{prop}

\begin{pf}
  By Theorem \ref{restr} and the structure of Artin-Schreier
  extensions, any $k$-isomorphism between these two curves
  must be of the type:
  \begin{equation}
    \label{eq:1}
    (x,y)\mapsto (\gamma(x),y+v(x)),
  \end{equation}
  for some $\gamma\in\pg,v(x)\in k(x)$, 
  and the identity $u'(x)=u(\gamma(x))+v(x)+v(x)^2$ must hold.\qed
\end{pf}

Therefore, the classification of hyperelliptic curves up to
$k$-isomor\-phism amounts to the classification of
rational functions $u(x)\in k(x)$
(with
sufficiently many poles of odd order) under the double action of
$\as(k(x))$ and $\pg$. In the next paragraph we shall indicate how to
carry out this classification, and in paragraph 1.3 we shall apply
this procedure for curves of genus two in a very explicit manner.

\subsection{Ramification divisor}

Given a smooth, projective curve $C$
defined over $k$ and a separable $k$-morphism $\pi\colon C\To
\prl$, we define the {\it ramification divisor} of $\pi$ as the
divisor, $\op{Diff}(C/\prl)$, of $C$ naturally associated to the
different of the extension of function fields $k(C)/\pi^*(k(\prl))$ 
(cf. \cite[ch.III]{Stich}). For any point $P\in C( \kb)$, if $e_P$ is the
ramification index at $P$ and $d_P$ is the exponent of the different
at $P$, then $d_P\ge e_P-1$, with equality if and only if $e_P$ is
odd (cf. \cite[III,5.1]{Stich}); hence, the support of the ramification
divisor is the set of ramification points of $\pi$.

We have seen in the last paragraph that any $u(x)\in k(x)-\as(k(x))$
determines a smooth projective curve $C_{u(x)}$ with quasi-affine
model $y^2+y=u(x)$ and implicitly equipped with a separable morphism of
degree 2 to $\prl$. By a procedure that goes back to Hasse
(cf. \cite{Hasse}), one can find a suitable element $v(x)\in k(x)$ such
that $u(x)+v(x)+v(x)^2$ has no poles of even order. Once this
normalization is achieved, the location and order of the poles of
$u(x)$ determine the ramification divisor of $C_{u(x)}$.
More precisely:

\begin{prop}\cite[III.7.8]{Stich}\label{ramdiv}
  Let $u(x)\in k(x)-\as(k(x))$ be a rational function having no poles
  of even order. Let $C=C_{u(x)}$ be the smooth, projective curve
  defined over $k$, with quasi-affine model $y^2+y=u(x)$. Then,
  $$
  \op{Diff}(C/\prl)= \sum_{Q\in \prl(\kb)} \Big(\sum_{P\mapsto Q}(m_Q+1
  ) P\Big),
  $$
  where
  $$
  m_Q = \left\lbrace
    \begin{array}{ll}
      -1,& \text{ if }\op{ord}_Q(u(x))\geq0,\\ m, & \text{ if
      }\op{ord}_Q(u(x))=-m<0.
    \end{array}
  \right.
  $$
  Moreover,
  $$k(C)\cap \kb =k\ \Leftrightarrow\ u(x)\not\in
  k+\as(k(x))\ \Leftrightarrow\ \op{Diff}(C/\prl)\ne 0,$$
  and, when this
  condition is satisfied, one has $\op{deg}(\op{Diff}(C/\prl))=2g+2$,
  where $g$ is the genus of $C$. In particular, if
  $\op{deg}(\op{Diff}(C/\prl))\ge 6$, then $C$ is a hyperelliptic curve
  defined over $k$.
\end{prop}

By Theorem \ref{restr}, if $C$ is a hyperelliptic curve, the
ramification divisor is independent of $\pi$ and, for any isomorphism
$\varphi\colon C\To C'$, we have
$\varphi^*(\op{Diff}(C'/\prl))=\op{Diff}(C/\prl)$. From now on we
shall denote:
$$
W:=\pi_{\ast}(\op{Diff}(C/\prl)),\quad \gw:=\{\ga \in
\pg \tq \ga^*(W)=W\}.
$$
These two objects depend on $\pi$, but the
reference to $\pi$ is omitted and it will be always implicit in the
context. Since any $k$-automorphism $\varphi$ of $C$ leaves the
ramification divisor invariant, the automorphism of $\prl$ associated
to $\varphi$ lies in $\gw$; hence, we have actually an exact
sequence:
$$
 1\To\io \To \op{Aut}_k(C) \stackrel{\pi}\To\gw.
$$

A natural strategy to classify hyperelliptic curves up to
$k$-isomor\-phism is to determine first, up to the action of $\pg$, the
possible divisors $W$ of $\prl$ that can appear as the push-forward
of ramification divisors of hyperelliptic curves and, afterwards,
classify the curves linked to a concrete divisor $W$. Proposition
\ref{ramdiv} takes care of the first step. Let us indicate a
procedure to carry out the second step.

Let $W$ be a fixed effective divisor of $\prl$ with even coefficients
and let $W'$ be the divisor obtained by lowering all positive
coefficients by one. Note that the isotropy subgroups of $W$, $W'$
under the action of $\pg$ are the same. Let ${\cal R}_{W'}\subseteq
k(x)$ be the set of rational functions having $W'$ as divisor of
poles and let $\n\subseteq{\cal R}_{W'}$ be a system of
representatives of these functions modulo $\as(k(x))$. We can define
an action of $\gw$ over $\n$ on the right. Given $u(x)\in \n$ and
$\ga\in\gw$, the divisor of poles of $\ga^*(u(x))=u(\ga(x))$ is again
$W'$ and we define $u^{\ga}(x)$ to be the only element of $\n$ such
that:
$$ u(\ga(x))\equiv u^{\ga}(x)\md{\as(k(x))}. $$

For any $\ga\in\gw$ and $u(x)\in\n$, denote by $v_{u,\ga}(x)\in k(x)$
any choice of a rational function satisfying: $$
u(\ga(x))=u^{\ga}(x)+v_{u,\ga}(x)+v_{u,\ga}(x)^2.$$

Let us check that $(\gamma,u(x))\mapsto u^{\ga}(x)$ defines an action
of $\gw$ on $\n$ indeed. 
For any couple $\ga,\eta\in\gw$, we have:
\begin{myeqnarray*}
(\ga\eta)^*(u(x))&=&u^{\ga\eta}(x)+v_{u,\ga\eta}(x)+v_{u,\ga\eta}(x)^2,\\
\eta^*(\ga^*(u(x)))&=&(u^{\ga})^{\eta}(x)+
v_{u^{\ga},\eta}(x)+v_{u^{\ga},\eta}(x)^2+\\
&&{}+\eta^*(v_{u,\ga}(x))+
\eta^*(v_{u,\ga}(x))^2.
\end{myeqnarray*}%
Since these two elements coincide, the functions $u^{\ga\eta}(x)$,
$(u^{\ga})^{\eta}(x)$ must be equal too, since both belong to $\n$
and they are congruent modulo $\as(k(x))$.

For any $u(x)\in\n$, denote by $\gu$ the isotropy group of $u(x)$:
$$\gu=\{\ga\in\gw \tq u(x)=u^{\ga}(x)\}.$$

\begin{prop}\label{pr}
  Let $W$, $\gw$, $\n\subseteq k(x)$ be as above.
  \begin{enumerate}
  \item For any $u(x),u'(x)\in\n$, the curves $C_{u(x)}$,
    $C_{u'(x)}$ are $k$-isomorphic if and only if there exists
    $\ga\in\gw$ such that $u'(x)=u^{\ga}(x)$.
  \item For any $u(x)\in\n$, we have an exact sequence:
    $$
    1\To\io \To \op{Aut}_k(C_{u(x)}) \stackrel{\pi}\To\gu\To 1.
    $$
    Moreover, if the map $v_u\colon\gu\To k(x)$, given by $\ga\mapsto
    v_{u,\ga}(x)$, is a homomorphism, then the exact sequence splits.
  \end{enumerate}
\end{prop}

\begin{pf}
  The first assertion is an immediate consequence of Proposition
  \ref{xy}. Let $\varphi$ be a $k$-automorphism of $C_{u(x)}$.
  By (\ref{eq:1}), the automorphism $\ga=\pi(\varphi)$ satisfies:
  $$u(x)\equiv u(\ga(x))\md{\as(k(x))}, $$
  so that $u(x)=u^{\ga}(x)$ and
  $\ga\in \gu$. Conversely, given $\ga\in\gu$, we can construct
  $\varphi$ by taking:
  $$ \varphi^*(x)=\ga(x),\qquad\varphi^*(y)=y+v_{u,\ga}(x). $$
  Clearly, if $v_u$ is a homomorphism,
  the map $\ga\mapsto\varphi$ obtained in this way is a homomorphism
  too.\qed
\end{pf}

\subsection{Normal equations for curves of genus two}
Let $C$ be a
smooth, projective, geometrically irreducible curve of genus two
defined over $k$. The canonical morphism, $\pi\colon C\To\prl$, is
defined over $k$ and has degree 2, so that $C$ is a hyperelliptic
curve. Hence, $C$ admits a quasi-affine model $y^2+y=u(x)$, for
certain $u(x)\in k(x)-\as(k(x))$ without poles of even order. By
Proposition \ref{ramdiv}, the ramification divisor of $C$ has degree
6 and one of the following possible forms:
\begin{myeqnarray*}
\lefteqn{\op{Diff}(C/\prl)=}\qquad\\&&
  =\left\lbrace \begin{array}{ll}
    2P_1+2P_2+2P_3, & \text{ if $u(x)$ has three simple poles},\\
    2P_1+4P_2,& \text{ if $u(x)$ has two poles, of orders 1,3},\\
    6P,& \text{ if $u(x)$ has only one pole, of order 5}.
  \end{array}\right.
\end{myeqnarray*}%
We label these possibilities respectively as: case (1,1,1), case
(1,3), and case (5). Since this divisor is defined over $k$, in cases
(1,3) and (5) all points in the support are defined over $k$ too.
However, in case (1,1,1) we have three different arithmetic
structures for the ramification divisor: three points defined over
$k$, one point defined over $k$ and two points conjugate over a
quadratic extension of $k$ or, finally, all three points conjugate
over a cubic extension of $k$. We shall refer to them respectively as
cases (1,1,1)-split, quadratic and cubic.

The divisor $W=\pi_{\ast}(\op{Diff}(C/\prl))$ has support on the
poles of $u(x)$.
In cases (1,1,1)-split, (1,3) and (5), by applying a suitable
$k$-automorphism of $\prl$, we can impose that this support 
is a subset of $\{\infty,0,1\}$. Let us display now a couple of
Lemmas that are necessary to fix the support of $W$ in the cases
(1,1,1)-quadratic and cubic.

\begin{lem}\label{le}
  Let $S,S'$ be two $G_{k}$-subsets of $\prl(\kb)$
  with $\abs{S}=\abs{S'}=3$. Then, there exists $\gamma \in \pg$
  such that $\gamma(S)=S'$ if and only if $S$ and $S'$ are isomorphic
  as $G_{k}$-sets.
\end{lem}

\begin{pf}
  Clearly, any $\gamma \in \pg$ respects the
  $G_k$-structure. Conversely, if $\rho\colon S \longrightarrow S'$ is
  a $G_{k}$-bijection, take $\gamma \in \op{PGL}_2(\kb)$ such that it
  coincides with 
  $\rho$ over $S$. Then, for any $\sigma \in G_k$ and any $x\in S$, we have:
  $$\ga (^{\sigma}\! x)=\rho (^{\sigma}\! x) =\,
  ^{\sigma}\!(\rho(x))=\,^{\sigma}\!(\ga(x))=
  \,^{\sigma}\!\ga(^{\sigma}\! x).
  $$
  Hence $\ga=\,^{\sigma}\!\ga$ for all $\sigma\in G_k$ and $\ga$ belongs
  to $\pg$.\qed
\end{pf}

\begin{lem}\label{cubic}
  For any cubic extension $K/k$, there exists a generator $\theta\in K$
  with minimal equation $x^3+sx+s$, $s\neq 0$, over $k$. If $w, w+1 \in
  \kb$ are the two roots of the equation $s+1=w+w^2$, then the other
  two roots of the polynomial $x^3+sx+s$ are:
  $$
  \theta'=\theta(\theta +w),\qquad  \theta''=\theta(\theta+w+1).
  $$
  In particular, the extension $K/k$ is cyclic if and only if
  $s+1\in\as(k)$.
\end{lem}

\begin{pf}
  Replacing $\theta$ by $\theta +\op{Tr}_{K/k}(\theta)$, we can
  assume that $\op{Tr}_{K/k}(\theta)=0$. Suppose that $x^3+ax+b$ is the
  minimal polynomial of $\theta$ over $k$. Replacing $\theta$ by
  $\theta+\theta^2$, if necessary, we can assume moreover that $a\ne
  0$. Then, the minimal polynomial of $ab^{-1}\theta$ is $x^3+sx+s$,
  where $s=a^3b^{-2}$. Finally, the assertions concerning the
  conjugates of $\theta$ can be checked by a direct computation.\qed
\end{pf}

We proceed now to exhibit families of normal equations describing all
$k$-isomorphism classes of curves of genus two. For each of the above
types of ramification divisor we fix a system of representatives of
divisors $W$ of this type modulo the action of $\pg$, and for each $W$
we apply the procedure indicated at the end of the last paragraph.
We find a suitable family $\n$ of normal equations 
and, just by
computing the action of $\gw$ on $\n$ and applying Proposition
\ref{pr}, 
we
classify the curves in $\n$ up to $k$-isomorphism and we determine the
structure of the $k$-automorphism group of each curve. In each case,
the rational functions in $\n$ depend on 
certain parameters, and we implicitly identify $\n$ with the set of
values of these parameters.

\subsubsection*{Case (1,1,1)-split}
Let $\n=k^*\times k^*\times k^*\times \left(\kas\right)$. Any curve
with ramification divisor of this type is $k$-isomorphic to one of
the curves $\CC$ with quasi-affine model:
$$ y^2 + y = ax+\frac bx+\frac c{x+1}+d ,\quad (a,b,c,d)\in\n. $$
For these curves,
$W=2[\infty ]+2[0]+2[1]$, and the set $Z_{\infty}$ of points at
infinity contains only the three Weierstrass points of the curve,
which are all defined over $k$. The subgroup $\g:=\gw$ is the
isotropy group of $\{\infty,\,0,\,1\}$ under the action of $\pg$ over
$\prl (\kb)$:
$$\g=\{x,\,1/x,\,1+x,\,x/(1+x),\,1/(1+x),\,(1+x)/x\}\simeq S_3.$$

In order to interpret properly the next result and the analogous
Propositions taking care of the other cases, we emphasize that, when
comparing two elements of $\n$, the fourth coordinate has to be
understood as an element of $\kas$.

\begin{prop}\label{uuu}
\begin{enumerate}
\item For any $(a,b,c,d),(a',b',c',d')\in \n$, the curves
  $\CCC$, $\CC$ are $k$-isomorphic if and only if $(a',b',c',d')$
  coincides with one of the following elements in the orbit of
  $(a,b,c,d)$ under the action of $\g$:
  \[
  (a,b,c,d)^{\ga} = 
  \left\lbrace \begin{array}{ll}
    (a,b,c,d),& \text{ if }\, \ga(x)=x, \\
    (b,a,c,d+c),& \text{ if }\, \ga(x)=1/x,\\
    (a,c,b,d+a),& \text{ if }\, \ga(x)=1+x,\\
    (c,b,a,d+c+b+a),& \text{ if }\, \ga(x)=x/(1+x),\\
    (b,c,a,d+c+b),& \text{ if }\, \ga(x)=1/(1+x),\\
    (c,a,b,d+b+a),& \text{ if }\, \ga(x)=(1+x)/x.
  \end{array}\right.
  \]
  
\item For any $(a,b,c,d)\in\n$, we have an split exact sequence:
  $$
  1\To\io\To \op{Aut}_k(\CC)\To \gab\To 1, $$
  $\gab$ being the
  isotropy group of $(a,b,c,d)$ under the action of $\g$:
  \[ \gab\simeq 
    \left\lbrace \begin{array}{ll}
      S_3,  & \text{ if }\, a=b=c \in \as(k),\\
      C_3,& \text{ if } \, a=b=c \notin \as(k), \\
      C_2, &\text{ if exactly two of the coefficients } a,b,c
      \text{ coincide}\\&\text{ and the third belongs to }\as(k), \\
      1,     & \text{ otherwise.}
    \end{array}\right.
  \]
  \end{enumerate}
\end{prop}

\begin{pf}
  To any $(a,b,c,d)\in\n$, we assign a unique element
  $u(x)=ax+\frac bx+\frac c{x+1}+d\in k(x)$, by making a choice of
  $d\in k$ inside the class modulo $\as(k)$ given by the fourth
  coordinate. The family of rational functions obtained in this way is
  a system of representatives modulo $\as(k(x))$ of the functions
  having  $[\infty]+[0]+[1]$ as divisor of poles and we can apply
  Proposition \ref{pr}. Only the fact that we can choose the mapping
  $v_u\colon \gu\To k(x)$ to be a homomorphism deserves some
  explanation. Given $u(x)$ as above and $\ga\in\gu$, it is easily
  checked that $v_{u,\ga}(x)$ is a constant satifying, respectively,
  $$v_{u,\ga}+v_{u,\ga}^2=0,c,a,c+b+a,c+b,b+a,$$
  according to the different possible values of $\ga$.
  Assume for instance
  that $a=b=c\in\as(k)$ and let $w\in k$ be such that $a=w+w^2$. Then,
  the following choice for $v_{u,\ga}$ determines a homomorphism:
  $$
  v_{u,\ga}= 
    \left\lbrace \begin{array}{ll}
      0,   & \text{ if }\, \ga\mbox{ has order 1 or 3}, \\
      w, & \text{ if }\, \ga\mbox{ has order 2}.
    \end{array}\right.
  $$
  The other cases are similar and easier.\qed
\end{pf}

From the exact sequence above, one can determine the $k$-automor\-phisms
of $\CC$.
For instance, if  $ a=b=c \in \as(k)$ and $a=w+w^2$ for certain $w\in
k$, the twelve $k$-automorphisms of $\CC$ are:
$$
\begin{array}{ll}
  (x,y) \stackrel{id} \longmapsto (x,y),&
  (x,y) \stackrel{\iota}  \longmapsto (x,y+1),\\
  (x,y) \stackrel{U} \longmapsto  (x+1,y+w), &
  (x,y) \stackrel{\iota U} \longmapsto (x+1,y+w+1),\\
  (x,y) \stackrel{V} \longmapsto (1/(1+x),y),&
  (x,y) \stackrel{\iota V} \longmapsto (1/(1+x),y+1),\\
  (x,y) \stackrel{V^2} \longmapsto ((1+x)/x,y),&
  (x,y) \stackrel{\iota V^2} \longmapsto ((1+x)/x,y+1),\\
  (x,y) \stackrel{VU} \longmapsto (1/x,y+w),&
  (x,y) \stackrel{\iota VU} \longmapsto (1/x,y+w+1),\\
  (x,y)\stackrel{V^2U} \longmapsto (x/(1+x),y+w), &
  (x,y)\stackrel{\iota V^2U} \longmapsto (x/(1+x),y+w+1).
\end{array}
$$

\subsubsection*{Case (1,1,1)-quadratic}
We fix a system of representatives $r\in k$ of all non-trivial
classes of $\kas$.

Let $\n=k^*\times \left((k\times k)-\{(0,\,0)\}\right)\times
\left(k/\as(k)\right)$. Any curve with ramification divisor of this
type is $k$-isomorphic to one of the curves $\CC^r$ with quasi-affine
model:
\begin{equation}\label{quad}
  y^2+y=ax+\frac{bx+c}{x^2+x+r}+d,
  \quad (a,b,c,d)\in\n.
\end{equation}
For these curves, $W=2[\infty]+2[\theta]+2[\theta']\,$, where $\theta,
\theta'\in\kb$ are the 
roots of the irreducible polynomial $x^2+x+r$. The set $Z_{\infty}$
of points at infinity contains only the three Weierstrass points of
the curve, one of them defined over $k$, the other two defined over
the quadratic extension determined by the class of $r$ modulo
$\as(k)$. The subgroup $\g:=\gw$ is the isotropy group of
$\{\infty,\,\theta,\,\theta'\}$ under the action of $\pg$ over $\prl
(\kb)$:
$$\g=\{x,\,1+x\}\simeq C_2.$$

By Lemma \ref{le}, two curves corresponding to different values of $r$
are not $k$-isomorphic. If $k$ is a finite field, there is a single
value of $r$: any choice of $r\in k-\as(k)$.

\begin{prop} Let $r\in k-\as(k)$ be fixed.
\begin{enumerate}
\item For any  $(a,b,c,d),(a',b',c',d')\in \n$, the curves
  $\CCC^r$, $\CC^r$ are $k$-isomorphic if and only if $(a',b',c',d')$
  coincides with one of the following elements in the orbit of
  $(a,b,c,d)$ under the action of $\g$:
  \[
  (a,b,c,d)^{\ga} = 
    \left\lbrace \begin{array}{ll}
      (a,b,c,d),& \text{ if }\, \ga(x)=x, \\
      (a,b,b+c,d+a),& \text{ if }\, \ga(x)=1+x.
    \end{array}\right.
  \]

\item For any $(a,b,c,d)\in\n$, we have an split exact sequence: $$
  1\To\io\To \op{Aut}_k(\CC^r)\To \gab\To 1, $$
  $\gab$ being the
  isotropy group of $(a,b,c,d)$ under the action of $\g$:
  \[ \gab\simeq 
    \left\lbrace \begin{array}{ll}
      C_2, &\text{ if } b=0, \ a\in\as(k),\\
      1,     & \text{ otherwise.}
    \end{array}\right.
  \]
\end{enumerate}
\end{prop}

\begin{pf}
  We can argue as in the proof of Proposition \ref{uuu}. When
  $b=0$ and $a=w+w^2$ for some $w\in k$, it is obvious that the choice
  $v_{u,\ga}=w$, for $\ga(x)=1+x$, determines a homomorphism from $\gu$
  to $k$.\qed
\end{pf}

\subsubsection*{Case (1,1,1)-cubic}
We fix a system of representatives, $\ck$, of cubic extensions of $k$
modulo isomorphism. By Lemma \ref{cubic}, we can identify $\ck$ with a
subset of $k$ containing representatives $s\in k$ of cubic
irreducible polynomials of the type $x^3+sx+s$.

Let $\n=\left((k\times k\times k)-\{(0,0,0)\}\right)\times
\left(\kas\right)$. Any curve with ramification divisor of this type
is $k$-isomorphic to one of the curves $\CC^s$ with quasi-affine
model:
$$y^2+y= \frac{ax^2+bx+c}{x^3+sx+s}+d,\quad (a,b,c,d)\in\n. $$
For these curves, $W=2[\theta]+2[\theta']+2[\theta'']$, where
$\theta, \theta', \theta''\in\kb$ are the roots of the irreducible
polynomial $x^3+sx+s$. The set $Z_{\infty}$ contains the three
Weierstrass points of the curve, which are defined over the splitting
field of this polynomial, together with two points,
$P_{\infty},P_{\infty}^{\iota}$, conjugate by the hyperelliptic
involution. These two points lie over the point with projective
coordinates $(x,y,z)=(1,0,0)$, which is an ordinary double point with
tangents $y^2+yz+dz^2=0$; thus, they are defined over $k$ if
$d\in\as(k)$, otherwise they are defined over the quadratic extension
determined by $d$.

By Lemma \ref{le} two curves corresponding to different values of $s$
are not $k$-isomorphic. If $k$ is a finite field, there is a single
value of $s$: any choice of $s\in k$ such that $x^3+sx+s$ is
irreducible.

If the cubic extension
$k(\theta)/k$ is not cyclic, then 
the isotropy group $\g:=\gw$ of the set $\{\theta,\theta',\theta''\}$
under the action of $\pg$ is trivial,
and we have as an immediate consequence
of Proposition \ref{pr}:

\begin{prop}\label{noncy}
Let $s\in \ck$ correspond to a non-cyclic cubic extension.
Then,
\begin{enumerate}
\item For any $(a,b,c,d),(a',b',c',d')\in \n$, the curves
$\CCC^s$, $\CC^s$ are $k$-isomorphic if and only if
$(a',b',c',d')=(a,b,c,d)$.
\item For any $(a,b,c,d)\in\n$, we have $\op{Aut}_k(\CC^s)=\io$.
\end{enumerate}
\end{prop}

In the cyclic case,
and by Lemma \ref{cubic},  we have:
\[ \g = \{x,\,\frac{wx+s}{x+w+1},\,\frac{(1+w)x+s}{x+w}\}\simeq C_3, \]
where  $s=1+w+w^2$  for some $w\in k$. In order to compute the action
of $\g$ we need the following observation, whose proof is
straightforward:

\begin{lem}\label{mz}
  For any $w\in k$, let us consider the matrix
  $$M(w)=
  \begin{pmat}
    w^2+w^3 & 1+w^3 & 1+w+w^5&w^2 \\
    1+w+w^2 & w^2 & 1+w^2+w^4&w \\
    1+w & w & w^3&1\\0&0&0&1
  \end{pmat}.
  $$
  Then, $M(w+1)=M(w)^2$ and $M(w)^3=I_3$. Moreover, if we let
  $M(w)$ operate on the right over $\n$ just by matrix multiplication:
  $$(a,b,c,d)^{M(w)}:=(a,b,c,d)M(w),$$
  then, the fixed points of this
  action are $(a,a,a(w+w^2),d)$, for arbitrary $a\in k^*$ and $d\in
  \kas$.
\end{lem}

\begin{prop}\label{u}
  Let $s\in k$ be such that $x^3+sx+s$ is irreducible and $s=1+w+w^2$
  for some $w\in k$. Then,
  \begin{enumerate}
  \item For any $(a,b,c,d),(a',b',c',d')\in \n$, the curves
    $\CCC^s$, $\CC^s$ are $k$-isomorphic if and only if $(a',b',c',d')$
    coincides with one of the following elements in the orbit of
    $(a,b,c,d)$ under the action of $\g$:
    \[
    (a,b,c,d)^{\ga} = 
      \left\lbrace \begin{array}{ll}
        (a,b,c,d),& \text{ if }\, \ga(x)=x, \\
        (a,b,c,d)M(w),& \text{ if }\, \ga(x)=\frac{wx+s}{x+w+1},\\
        (a,b,c,d)M(w+1),& \text{ if }\, \ga(x)=\frac{(1+w)x+s}{x+w}.
      \end{array}\right.
    \]
  \item For any $(a,b,c,d)\in\n$, we have an split exact sequence:
    \begin{equation}\label{aut}
      1\To\io\To \op{Aut}_k(\CC^s)\To \gab\To 1,
    \end{equation}
    $\gab$
    being the isotropy group of $(a,b,c,d)$ under the action of $\g$:
    $$
    \gab\simeq 
    \left\lbrace \begin{array}{ll}
      C_3, &\text{ if } a=b, \ c=a(1+s),\\
      1,     & \text{ otherwise.}
    \end{array}\right.
    $$
  \end{enumerate}
\end{prop}

\begin{pf}
  We can argue as in the proof of Proposition \ref{uuu}. The
  computation of $\gab$ is consequence of Lemma \ref{mz}. When $a=b$,
  $c=a(1+s)$ we can choose $v_u\colon \gu\To k$ to be the trivial
  homomorphism.\qed
\end{pf}

For instance, when $ a=b,\,c=a(1+s)$, the  subgroup of
$\op{Aut}_k(\CC^s)$ generated by $(x,y) \mapsto
(\frac{wx+s}{x+w+1},y)$ is cyclic of order three.

If a cubic extension admits a generator with minimal polynomial over
$k$ of the type $x^3+s$, it might be easier to work with models
$\CC^s$ of the type:
$$y^2+y=\frac{ax^2+bx+c}{x^3+s}+d,\quad
(a,b,c,d)\in \n. $$

If the cubic extension is not cyclic, Proposition \ref{noncy} is
valid for these models as well. If a 
cyclic cubic field admits a cyclic generator, then $\mathbb F_4\subseteq k$.
In this case,
Proposition \ref{u} might be substituted by:

\begin{prop}
  Assume that $\mathbb F_4\subseteq k$ and let $s\in k$ be such that
  $x^3+s$ is irreducible. Then,
  \begin{enumerate}
  \item For any $(a,b,c,d),(a',b',c',d')\in \n$, the curves
    $\CCC^s$, $\CC^s$ are $k$-isomorphic if and only if
    $$ (a',b',c',d')=(\epsilon^2 a,\epsilon b,c,d)$$
    for some $\epsilon\in \mu_3$.
  \item For any $(a,b,c,d)\in\n$, we have an split exact sequence
    (\ref{aut}) where,
    \[ \gab\simeq \left \lbrace
      \begin{array}{ll}
        C_3, &\text{ if } a=b=0,\\
        1,     & \text{ otherwise.}
      \end{array}
    \right. \]
  \end{enumerate}
\end{prop}

\subsubsection*{Case (1,3)}
Let $\n=k^*\times k\times k^*\times (k/\as(k))$.
 Any curve with ramification divisor of this type is
$k$-isomorphic to one of the curves $\CC$ with quasi-affine model:
$$y^2+y=ax^3+bx+\frac{c}{x} + d, \quad (a,b,c,d)\in\n.$$

For these curves, $W=4[\infty]+2[0]$, and $Z_{\infty}$ contains only
the two Weierstrass points of the curve, which are defined over $k$.
The subgroup $\g:=\gw$ is, in this case, $\g=\{\la x\tq\la\in k^* \}
\simeq k^*$.

\begin{prop}
\begin{enumerate}
\item For any $(a,b,c,d),(a',b',c',d')\in \n$, the curves
  $\CCC$, $\CC$ are $k$-isomorphic if and only if $(a',b',c',d')$
  coincides with one of the following elements in the orbit of
  $(a,b,c,d)$ under the action of $\g$:
  $$ (a,b,c,d)^{\la x} = (\la^3
  a,\la b,\la^{-1}c,d),\quad \la\in k^*.$$
\item For any $(a,b,c,d)\in\n$, we have $\op{Aut}_k(\CC)=\io$.
\end{enumerate}
\end{prop}

\subsubsection*{Case (5)}
This case is known as the {\it supersingular} case, since the
jacobian of these curves splits up to isogeny as the square of a supersingular
elliptic curve over $\kb$. This case has been thoroughly studied by van der
Geer and van der Vlugt in the
finite field case (cf. \cite{VanderGeer2,VanderGeer}).

Let $\n=k^*\times k \times k \times (k/\as(k))$. Any curve with
ramification divisor of this type is $k$-isomorphic to one of the
curves $\CC$ with quasi-affine model:
\begin{equation}\label{ss}
  y^2 +y = a x^5 + b x^4 + c x^3 +  d, \quad (a,b,c,d)\in\n.
\end{equation}

For these curves, $W=6[\infty]$, and $Z_{\infty}$ contains only the
Weierstrass point of the curve, which is defined over $k$. The
subgroup $\g:=\gw$ is now the affine group,
$$\g=\{ \la x + \nu \tq (\la,\nu) \in k^*\times k \}\simeq k^* \rtimes k.$$

\begin{lem}\label{delta}
Given $a,b,c\in k$, $a\ne 0$, consider the linear separable
polynomial $E_{ac}(x)=a^4x^{16}+c^4x^8+c^2x^2+ax\in k[x]$. Then, the
subset
$$\Gamma_{abc}:=\left\lbrace\begin{array}{ll}
    \{(1,\nu)\in k^*\rtimes k\tq E_{ac}(\nu)=0\},&\mbox{ if } c\ne 0,\\
    \{(\la,\nu)\in k^*\rtimes k\tq \la^5=1,\ E_{ac}(\nu)=b(1+\la)\},&
    \mbox{ if } c=0,
 \end{array}\right.$$
is a subgroup of $k^*\rtimes k$ and the map
$$
\delta_{abc}\colon\Gamma_{abc}\To \kas, \quad (\la,\nu)\longmapsto
a\nu^5+b\nu^4+c\nu^3,
$$
is a group homomorphism.
\end{lem}

\begin{pf} Given any $(\la,\nu),(\la',\nu')\in \Gamma_{abc}\subseteq
k^*\rtimes k$, their product is $(\la\la',\la\nu'+\nu)$. If $c\ne 0
$, $\Gamma_{abc}$ is a subgroup by the linearity of $E_{ac}(x)$. If
$c=0$ we have also:
$$
E_{ac}(\la\nu'+\nu)=\la E_{ac}(\nu')+E_{ac}(\nu)=
\la b(1+\la')+b(1+\la)=b(1+\la\la').
$$
The assertion concerning $\delta_{abc}$ is left to the reader.\qed\end{pf}

We note that if $c\ne 0$ or $\mu_5(k)=\{1\}$ then
$\op{Ker}(\delta_{abc})$ is isomorphic to a subgroup of $C_2\times
C_2\times C_2\times C_2$.

As in the previous cases, Proposition \ref{pr} leads to:

\begin{prop}\label{cas5}
\begin{enumerate}
\item For any $(a,b,c,d),(a',b',c',d')\in \n$, the curves
$\CCC$, $\CC$ are $k$-isomorphic if and only if $(a',b',c',d')$
coincides with one of the following elements in the orbit of
$(a,b,c,d)$ under the action of $\g$:
\[
(a,b,c,d)^{\la x+\nu} =(\la^5 a, \la^4(b+E_{ac}(\nu)),\la^3
c,a\nu^5+b\nu^4+c\nu^3+d).
\]
\item 
  The
  isotropy group of any $(a,b,c,d)\in\n$ under the action of $\g$
  is $\op{Ker}(\delta_{abc})$.
  We have an
  exact sequence (which rarely splits):
\begin{equation}\label{meta}
 1\To\io\To \op{Aut}_k(\CC)\To
\op{Ker}(\delta_{abc})\To 1.
\end{equation}
\end{enumerate}
\end{prop}

\subsection{Number of curves of genus two over a finite field}
Let
$k=\ff{}$ be a finite field with $q=2^m$ elements. In this paragraph
we find an explicit formula for the number of curves of genus two
defined over $k$ up to $k$-isomorphism. Actually, we count how many
curves there are for each type of ramification divisor. After the
results of paragraph 1.3, the number of curves having a concrete type
of ramification divisor coincides with the numbers of orbits of the
action of certain finite group $\g:=\gw$ on certain finite set $\n$.
Denoting by $\g \backslash \! \backslash \n $ the set of
these orbits, we can count its cardinality by a well-known formula:

\begin{equation}\label{comb}
  \vert \g \backslash \! \backslash \n \vert =
  \frac{1}{\vert \g \vert} \sum_{\ga\in \g} \vert \n_{\ga}\vert =
  \sum_{\ga\in {\cal C}} \frac{\vert \n_{\ga}\vert}{\vert
    \g_{\ga}\vert},
\end{equation}
where ${\cal C}$ is a system of
representatives of all conjugation classes of $\g$ and
$$
\n_{\ga}=\{x\in \n \tq x^{\ga}=x\}, \quad \g_{\ga}= \{\rho \in \g
\tq \rho \ga = \ga \rho \}.
$$

From now on, a summand $[a]_{b|m}$ in a formula means ``add $a$ if $b$
divides $m$''. 

\renewcommand\arraystretch{1.4}
\begin{thm}
  There are $2q^3+q^2+q-2+[8]_{4|m}$ $k$-isomorphism classes of curves
  of genus two defined over $\ff{}$. The number of curves with a given
  type of ramification divisor is given in the following table:
  $$\begin{array}{lc}
      \hline
      \text{(1,1,1)-split}
      &
      \frac 16q(q-1)(2q-1)
      \\
      \text{(1,1,1)-quadratic}
      &
      \frac 12 (q-1)(2q^2+q-4)
      \\
      \text{(1,1,1)-cubic}
      &
      \frac 13(2q^3+4q-6)
      \\
      \text{(1,3)}
      &
      2q(q-1)
      \\
      \text{(5)}
      &
      4q-2+[8]_{4|m}
      \\\hline
    \end{array}
    $$\phantom{x}
\end{thm}

\begin{pf}
  By (\ref{comb}), we need only to compute $\vert\n_{\ga}\vert$,
  $\vert\g_{\ga}\vert$ in all cases.
  
  In the case (1,1,1)-split, $\g\simeq S_3$ and we can choose ${\cal
    C}=\{1,\tau,\sigma\}$, where $\tau(x)=1+x$, $\sigma(x)=1/(1+x)$.
  Clearly, $\g_{\tau}=\{1,\tau\}$, $\g_{\sigma}=\{1,\sigma,\sigma^2\}$,
  and
  $$\n_{\tau}=\{(a,b,b,d)\tq a\in\as(k)-\{0\}, b\in k^*, d\in
  \kas\},$$
  $$\n_{\sigma}=\{(a,a,a,d)\tq a\in k^*,d\in \kas\}. $$
  Hence,
  $\vert \n_1 \vert=\vert \n \vert=2(q-1)^3$, $\vert \n_{\tau} \vert=
  2(q-1)(q/2 -1)$ and $\vert \n_{\sigma} \vert=2(q-1)$.
  
  In the case (1,1,1)-quadratic, $\g=\{1,\tau\}$, where $\tau(x)=1+x$;
  hence,
  $$\n_{\tau}=\{(a,0,c,d)\tq a\in\as(k)-\{0\}, c\in k^*, d\in
  \kas\}.
  $$
  Therefore, $\vert \n_1 \vert=2(q-1)(q^2-1)$ and
  $\vert \n_{\tau} \vert=2(q-1)(q/2-1)$.

  In the case (1,1,1)-cubic, $\g=\{1,\ga,\ga^2\}\simeq C_3$ and
  $$\n_{\ga}=\n_{\ga^2}=\{(a,a,a(1+s),d)\tq a\in k^*, d\in \kas\}.$$
  Hence, $\vert \n_1 \vert=2(q^3-1)$ and $\vert \n_{\ga}
  \vert=\vert \n_{\ga^2} \vert=2(q-1)$.

  In the case (1,3), $\g\simeq k^*$ and for any $\la\in k^*$, $\la\ne
  1$, we have $\n_{\la}=\emptyset$. Therefore, only
  $\vert \n_1 \vert=2q(q-1)^2$ counts.

  In the case (5), $\g= k^*\rtimes k$ and we can choose: ${\cal
    C}=\{\tau\}\cup\{\sigma_{\la}\tq \la\in k^*\}$, where $\tau(x)=1+x$,
  $\sigma_{\la}(x)=\la x$. We have clearly (if $\la\ne 1$),
  $$\Gamma_{\tau}=\{x+\nu\tq \nu\in k\},\quad
  \Gamma_{\sigma_{\la}}=\{\la'x\tq\la'\in k^*\}.$$
  
  For any given $a\ne 0$, there are only two values of $c$ such that
  $E_{ac}(1)=0$; these are: $c_0=\sqrt{a+a^2}$, $c_1=1+\sqrt{a+a^2}$.
  Hence,
  $$ \n_{\tau}=\{(a,b,c,d)\tq a\in k^*,c\in\{c_0,c_1\},b\in
  a+c+\as(k),d\in \kas\}, $$
  and $\vert\n_{\tau}\vert=2q(q-1)$.
  Finally, if $\la\not\in \mu_5(k)$, we have
  $\n_{\sigma_{\la}}=\emptyset$, whereas for $\la\in \mu_5(k)$, $\la\ne
  1$, we have $\n_{\sigma_{\la}}=\{(a,0,0,d)\tq a\in k^*, d\in \kas\}$,
  with cardinality $2(q-1)$, whereas $\vert \n_1 \vert=2(q-1)q^2$.\qed
\end{pf}

Also, from the results of paragraph 1.3 we can obtain a mass formula
for the number of curves of genus two:

\begin{thm}\label{mass}
  Let $[C]$ run on the $k$-isomorphism classes of
  smooth projective curves $C$ of genus two defined over $k$. Then,
  $$\sum_{[C]} \vert \op{Aut}_k(C)\vert^{-1}=q^3.$$
  More precisely, the
  partial weighted sums $\sum_{[C]} \vert \op{Aut}_k(C)\vert^{-1}$,
  letting $[C]$ run on the $k$-isomorphism classes of curves having a
  fixed type of ramification divisor, are given in the following
  table:
  $$\begin{array}{lc}
  \hline
  \text{(1,1,1)-split}
  &
  \frac 16(q-1)^3
  \\
  \text{(1,1,1)-quadratic}
  &
  \frac 12(q-1)(q^2-1)
  \\ 
  \text{(1,1,1)-cubic}
  &
  \frac 13(q^3-1)
  \\
  \text{(1,3)}
  &
  q^2-q
  \\
  \text{(5)}
  &
  q
  \\\hline
\end{array}
$$\phantom{x}
\end{thm}

\begin{pf}
  Each partial sum can be computed as
  \begin{myeqnarray*}
  \lefteqn{\sum_{(a,b,c,d)\in \g
    \backslash \! \backslash \n} \frac 1{\vert
    \op{Aut}_k(\CC)\vert}=}\qquad\\&&=\sum_{(a,b,c,d)\in \g \backslash \!
    \backslash \n} \frac 1{2\vert \gab\vert}= \sum_{(a,b,c,d)\in \n}
  \frac 1{2\vert \g\vert}=\frac {\vert \n\vert}{2\vert\g\vert},
\end{myeqnarray*}%
  for
  certain finite set $\n$ and finite group $\g$ acting on $\n$,
  specified in paragraph 1.3. \qed
\end{pf}

Both the total and the weighted number of supersingular curves had
been already obtained by van der Geer and van der Vlugt in
\cite{VanderGeer}.

\renewcommand\arraystretch{1.}

\section{Geometric invariants of curves of genus two}
In this section we recall from \cite{Igusa} the definition of
invariants that classify curves of genus two over $\kb$ up to
isomorphism and we compute explicitly these invariants in terms of
the rational models introduced in section 1. In paragraph 2.2 we deal
with rationality questions concerning the curve $C$, its isomorphism
class and the values of their geometric invariants.

A $\kb$-isomorphism between two curves defined over $k$ is simply
called an {\it isomorphism} and we denote by $\au:=\op{Aut}_{\kb}(C)$
the full group of automorphisms of a curve $C$.

\subsection{Igusa invariants}
Let $C$ be a curve of genus two defined over $\kb$.
According to the different possibilities for the ramification divisor
of $C$, the curve admits a quasi-affine model of the type:
$$
\begin{array}{ll}
  y^2+y=ax+\displaystyle\frac bx+\displaystyle\frac c{x+1},\quad abc\ne 0,&\qquad\text{type (1,1,1),}\\[1em]
  y^2+y=ax^3+bx+\displaystyle\frac cx,\quad ac\ne 0,&\qquad\text{type (1,3),}\\[1em]
  y^2+y=ax^5+cx^3,\quad a\ne 0,&\qquad\text{type (5),}
\end{array}%
$$
and we define, respectively,
$$
\begin{array}{l@{\qquad}l@{\qquad}l}
  j_1(C):=abc,& j_2(C):=ab+bc+ca, & j_3(C):=a+b+c, \\[1em]
  j_1(C):=0,& j_2(C):=ac^3, & j_3(C):=bc,\\[1em]
  j_1(C):=0,& j_2(C):=0, & j_3(C):=c^5/a^3.
\end{array}
$$



\begin{defn}
For any curve $C$ of genus two defined over $\kb$ we define the
$j$-invariant of $C$ as:
$$ j(C):=(j_1(C),j_2(C),j_3(C))\in \kb^3. $$
\end{defn}

By applying the results of paragraph 1.3 to the field $\kb$ we get:

\begin{thm}[Igusa]\label{igusa}\noindent
\begin{enumerate}
\item The invariant $j(C)$ depends only on the isomorphism class of $C$.
\item Two curves of genus two defined over $\kb$ are isomorphic if and
  only if they have the same $j$-invariant.
\item The following table sums up the ramification type of $C$ and the
possible structures of $\au$ in terms of the invariant
$j(C)=(j_1,j_2,j_3)$:
$$
\begin{array}{llll}
  \text{Type} & \text{Condition} & \au &
  \text{Condition}\\\hline
  (1,1,1) & j_1\neq 0 & C_2 & j_1\neq j_2j_3 \\ 
  & & C_2\times C_2 & j_1=j_2j_3,\ j_1\neq j_3^3 \\ 
  & & C_2\times S_3 & j_1=j_2j_3,\ j_1=j_3^3 \\ \hline
  (1,3) & j_1=0,\ j_2\neq0 & C_2 & - \\ \hline
  (5) & j_1=j_2=0 & M_{32} & j_3\neq 0 \\ 
  & & M_{160} & j_3=0 \\ \hline
\end{array}
$$
\end{enumerate}
\end{thm}

The groups $M_{32}$, $M_{160}$ sit in the middle of the non-split
exact sequence 
(\ref{meta})
and we shall study their structure more
closely in section 3.
Over an algebraically closed field, the group
$\op{Ker}(\delta_{abc})$ of (\ref{meta}) is clearly:
$$
\op{Ker}(\delta_{abc})=\left\lbrace
\begin{array}{ll} \op{Ker}E_{ac}\simeq C_2\times C_2\times C_2\times
 C_2,&\mbox{ if }c\ne 0\\\mu_5\rtimes \op{Ker}E_{ac},&\mbox{ if }c=0.
\end{array} \right.
$$

For $j_1\ne 0$, the first three conditions of the last column of the
above table are respectively equivalent to the fact that the
polynomial $x^3+j_3x^2+j_2x+j_1\in \kb[x]$ is separable, has a double
root or a triple root.

We compute now, in terms of the rational models given in section 1,
the $j$-invariant of a curve of genus two defined over 
$k$.

In all
cases, the curve determined by the parameters $(a,b,c,d)$ is
isomorphic to the curve determined by $(a,b,c,0)$; hence, the fourth
parameter $d$ may be ignored in the computation of the $j$-invariant.
Moreover, in the supersingular case, the curve (\ref{ss}) determined
by the parameters $(a,b,c,d)$ is isomorphic to the curve determined
by $(a,0,c,0)$. After these remarks, it only remains to treat the
cases (1,1,1)-quadratic or cubic.

\begin{lem}\label{qc}
\begin{enumerate}
\item For a curve $C$ of type (1,1,1)-quadratic with quasi-affine model
given by equation (\ref{quad}), we have:
$$
j_1(C)=a(rb^2+bc+c^2),\quad j_2(C)=ab+rb^2+bc+c^2,\quad j_3(C)=a+b.
$$
\item Let $s,t\in k$ be such that the polynomial $x^3+tx+s$ is irreducible in
$k[x]$. For a curve $C$ of type (1,1,1)-cubic with quasi-affine
model
\begin{equation}\label{rara}
  y^2+y=\frac{ax^2+bx+c}{x^3+tx+s},
\end{equation}
the $j$-invariant of $C$ is given by:
\begin{myeqnarray*}
j_1(C)&=&\frac 1{s^3}(a^2(s^2a+stb+t^2c)+bc(tb+sa)+sb^3+c^3),\\
j_2(C)&=&\frac 1{s^4}((a^2t^2+abs+c^2)(s^2+t^3)+cst(tb+sa)+b^2t^4),\\
j_3(C)&=&\frac 1s(ta+c).
\end{myeqnarray*}%
\end{enumerate}
\end{lem}

\begin{pf}
  Assume that $C$ is given by equation (\ref{quad}). If
  $\theta,\theta'$ are the roots of $x^2+x+r$, we have that
  $$
  \frac{bx+c}{x^2+x+r}=
  \frac{b\theta+c}{x+\theta}+\frac{b\theta'+c}{x+\theta'},
  $$
  hence, the mapping $(x,y)\mapsto (x+\theta,y)$ is an isomorphism
  from $C$ to the curve $C'$ given by an split equation
  $y^2+y=a'+(b'/x)+(c'/(x+1))+d'$, where
  $$ a'=a,\quad b'=b\theta+c,\quad c'=b\theta'+c. $$
  It is easily checked now that
  $j(C)=j(C')$ has the claimed values.

  Assume now that $C$ is given by equation (\ref{rara}). If
  $\theta,\theta',\theta''$ are the roots of $x^3+tx+s$, we have that
  $$
  \frac{ax^2+bx+c}{x^3+tx+s}=\frac
  A{x+\theta}+\frac{A'}{x+\theta'}+\frac {A''}{x+\theta''},
  $$
  where
  $A,A',A''\in\overline{k}$ are given by
  $$ \begin{pmat1} A\\A'\\A'' \end{pmat1}=
  \frac 1s \begin{pmat3}
    \theta^3&\theta^2&\theta\\
    {\theta'}^3&{\theta'}^2&\theta'\\
    {\theta''}^3&{\theta''}^2&{\theta''}
  \end{pmat3}
  \begin{pmat1} a\\b\\c \end{pmat1}.
  $$
  Hence, the
  mapping $(x,y)\mapsto
  (\frac{\theta'}{\theta}\frac{x+\theta'}{x+\theta},y)$ is an
  isomorphism from $C$ to the curve $C'$ given by the  equation
  $y^2+y=a'+(b'/x)+(c'/(x+1))+d'$, where
  $$ (a',b',c')=\frac
  1s\left(\theta^2A,(\theta')^2A',(\theta'')^2A''\right).
  $$
  Therefore,
  the values of $j_i(C)=j_i(C')$, $i=1,2,3$, are the coefficients of
  the minimal polynomial of $(a \theta^5+b \theta^4+c\theta^3)/s^2$
  over $k$.\qed
\end{pf}

\subsection{Field of moduli and field of definition}
Let $C$ be a curve of genus two defined over $\kb$. An intermediate
field $K$ of $\kb/k$ is said to be a {\it field of moduli} for $C$ if
$C$ is isomorphic to $^{\sigma}\! C$ for all $\sigma\in
G_K$. By Theorem \ref{igusa},
$$C\simeq\, ^{\sigma}\!
C\Longleftrightarrow j(C)=j(^{\sigma}\!
C)\Longleftrightarrow j(C)=\, ^{\sigma}\! j(C),$$
so that
$K$ is a field of moduli for $C$ if and only if  $j(C)\in K^3$.

If the curve $C$ is defined over $k$, then clearly $j(C)\in k^3$ and
$k$ is a field of moduli for $C$. We want to see the converse: given
a curve of genus two defined over $\kb$, having $k$ as a field of
moduli, then $C$ is isomorphic to a curve defined over $k$. This
result is false for curves of genus two over fields of characteristic
different from 2, with non-trivial two-torsion in the Brauer group,
but the situation is well understood (cf. \cite{Mestre}, \cite{CaQ}).

\begin{thm}\label{cosdef}
For any given value of $j\in k^3$, there exists a curve of genus two
$C$ defined over $k$ such that $j(C)=j$.
\end{thm}

Equivalently, we could formulate this result in the following form:

\begin{cor}\label{jc}
Let ${\cal C}$ denote the quotient set of all curves of genus two
defined over $k$, classified up to $\kb$-isomorphism. Then, the
$j$-invari\-ant sets a bijection between ${\cal C}$ and $k^3$.
\end{cor}

For $k$ a finite field, Theorem \ref{cosdef} follows immediately from
Theorem \ref{mass} and the following result of van der Geer and van der
Vlugt
(cf. \cite[Proposition 5.1]{VanderGeer}):
$$\sum_{[C]}\vert\op{Aut}_k(C)\vert^{-1}=\vert
j({\cal C})\vert,$$
where $[C]$ runs over all $k$-isomorphism classes
of curves of genus two defined over $k$. We give now a direct proof
of Theorem \ref{cosdef}, for general $k$, essentially by exhibiting
rational models of curves of genus two with prescribed $j$-invariant.

\begin{pf*}{Proof of Theorem \ref{cosdef}}
  Let $j=(j_1,j_2,j_3)\in k^3$ be
  given. We find a quasi-affine model for a curve $C$ defined over $k$
  such that
  $j(C)=j$.

  If $j_1=j_2=0$, the following supersingular curves have $j(C)=j$:
  $$
    \begin{array}{ll}
      y^2+y=\sqrt{j_3}\,x^5+\sqrt{j_3}\,x^3,&\mbox{ if }j_3\ne 0,\\
      y^2+y=x^5,&\mbox{ if }j_3=0.
    \end{array} 
  $$

  If $j_1=0$, $j_2\ne 0$, then the following curve of type (1,3)
  has $j(C)=j$: 
  $$ y^2+y=j_2x^3+j_3 x+\frac 1x. $$
  
  If $j_1\ne0$, consider the polynomial $F(x)=x^3+j_3x^2+j_2x+j_1\in
  k[x]$. We proceed in this case according to the arithmetic structure
  of the set of roots of $F(x)$. If $F(x)$ splits in $k[x]$,
  $F(x)=(x+a)(x+b)(x+c)$, then the curve $C$ of type (1,1,1)-split with
  quasi-affine model:
  $$ y^2+y=ax+\frac bx+\frac c{x+1}, $$
  has $j(C)=j$. If $F(x)=(x+a)Q(x)$ for some quadratic irreducible
  polynomial $Q(x)=x^2+ux+v$, then the curve $C$ of type
  (1,1,1)-quadratic with quasi-affine model:
  $$ y^2+y=ax+\frac{ux+u}{x^2+x+(v/u^2)}, $$
  has $j(C)=j$, by Lemma
  \ref{qc}. Finally, assume that $F(x)$ is irreducible in $k[x]$; let
  $\omega\in \overline{k}$ be a root of $F(x)$ and $K=k(\omega)$. The
  element $\theta=\omega+j_3$ has minimal polynomial $x^3+tx+s$ over
  $k$, where $t=j_2+j_3^2$ and $s=j_1+j_2j_3$. It is easy to check that
  $\theta^3,\theta^4,\theta^5$ is a $k$-basis of $K$; hence, there
  exist unique $a,b,c\in k$ such that
  \begin{equation}\label{fi}
    a\theta^5+b\theta^4+c\theta^3=s^2\omega.
  \end{equation}
  For these values of $a,b,c$, 
  the curve $C$ of type (1,1,1)-cubic with
  quasi-affine model given by (\ref{rara}) has $j(C)=j$,
  by the proof of Lemma \ref{qc}.
For the sake of completeness, we 
give an explicit equation (\ref{rara}) in terms of
$j_1,j_2,j_3$ in this last case.
If $j_2=j_3^2$, we can take $t=0$,
$s=j_1+j_2j_3$ and $$a=0,\quad b=j_1+j_2j_3, \quad
c=j_3(j_1+j_2j_3).$$ If $j_2\ne j_3^2$, taking
$t=s=\frac{(j_2+j_3^2)^3}{(j_1+j_2j_3)^2},$
the element
$\theta=\frac{j_2+j_3^2}{j_1+j_2j_3}(\omega+j_3)$ has minimal
polynomial $x^3+sx+s$ over $k$ (see the proof of Lemma \ref{cubic}).
For this choice of $t,s,\theta$, the values:
\begin{myeqnarray*}
a=b&=&\frac{(j_1+j_3^3)(j_2+j_3^2)^2}{(j_1+j_2j_3)^2},\\
  c&=&\frac{(j_2+j_3^2)^3(j_1(j_2+j_3^2)^2+j_3(j_1+j_3^3)^2)}{(j_1+j_2j_3)^4},
\end{myeqnarray*}%
satisfy (\ref{fi}).
\qed
\end{pf*}

\section{Arithmetic invariants of curves of genus two}
Let $C$ be a curve of genus two defined over $k$. A twist of $C$ over
$k$ is any curve of genus two, defined over $k$, which is isomorphic
to $C$. We denote by $\tw$ the quotient set of all twists of $C$ over
$k$, classified up to $k$-isomorphism. There is a well-known
isomorphism of pointed sets,
\begin{equation}
  \label{tw}\tw\stackrel{\sim}\To \h{\au},
\end{equation}
sending any twist $\phi\colon C'\stackrel{\sim}\To C$ to the
1-cocycle $\{\phi\,^{\sigma}\!\phi^{-1}\}_{\sigma\in
G_k}$. For details concerning non-abelian cohomology we address the
reader to \cite[Ch.VII, Annexe]{Se}.

In this section we obtain an explicit description of the set $\tw$ of
twists of a given curve $C$, in terms of the rational models of
section 1. This is achieved by an explicit computation of $\h{\au}$
for any possible structure of $\au$. The parameters describing
the pointed set $\h{\au}$ are
what we call the {\it arithmetic invariants} of
the curve. In this way, we obtain an explicit parameterization of all
$k$-isomorphism classes of curves of genus two, each class being
determined by a couple of invariants, one geometric and the other
arithmetic. For $k$ a finite field, this provides
formulas for the number of
curves whose full group of automorphisms has a concrete structure.
For finite fields of odd characteristic an analogous result has been
obtained in \cite{Ca}.

As a by-product, we are able to
implement an algorithm that builds up, almost directly, a faithful
and complete system of representatives of $k$-isomorphism classes of
curves of genus two over a finite field $k$ of even characteristic.
This algorithm was used in \cite{MN} to carry out a numerical
exploration about the existence of jacobians of curves of genus two
in a certain isogeny class of abelian surfaces over $k$.

\subsection{Hyperelliptic twists}
Let $C$ be a curve of genus two defined over $k$, with hyperelliptic
involution $\iota$. By Artin-Schreier theory, we have an isomorphism:
\begin{equation}\label{hac}
H:=\kas\simeq \h{\io}=\op{Hom}(G_k,\io).
\end{equation}
Thus,
elements in $H$ can be thought as 1-cocycles of $G_k$ with values in
$\au$ and they furnish twists of $C$ by (\ref{tw}). We obtain in this
way an action of the group $H$ on the set of curves of genus two
defined over $k$. The curves in the orbit of $C$ under this action
are called {\it hyperelliptic twists} of $C$, and we denote by
$\twc$ the quotient set of this orbit modulo $k$-isomorphism.

If $C$ admits a quasi-affine model $y^2+y=u(x)$ for certain
$u(x)\in k(x)$, the hyperelliptic twist of $C$ corresponding to
$d\in H$ is the curve $C'$ with model $y^2+y=u(x)+d$. If $v\in\kb$
satisfies $v^2+v=d$, the mapping $(x,y)\mapsto (x,y+v)$ is an
isomorphism between $C'$ and $C$, defined over the quadratic (or
trivial) extension of $k$ determined by $d$.

Any isomorphism, $\phi\colon C'\stackrel{\sim}\To C$, between two
curves of genus two commutes with the respective hyperelliptic
involutions: $\iota \phi=\phi\iota'$. Therefore, the action of $H$ is
compatible with isomorphisms and with $k$-isomorphisms. In particular,
we have an induced action of $H$ on the set $\tw$. For any $C'\in\tw$
we shall denote by $H_{C'}$ the isotropy group of $C'$ under this
action; note that the orbit of $C'$ can be identified with
$\op{Tw}^0(C'/k)$. In order to describe $\tw$ we need to compute the
set of orbits and all the isotropy subgroups.

Consider the central exact sequence of $G_k$-groups:
\begin{equation}\label{a} 1\To\io\To A\To A'\To 1,
\end{equation}
where $A:=\au$, $A':=A/\io$ denote respectively the full group and
the reduced group of automorphisms of $C$. Since $H^2(G_k,C_2)=0$,
the cohomology exact sequence of (\ref{a}) induces an exact
sequence of pointed sets,
\begin{equation}\label{b}
1\To\h{\io}/\delta((A')^{G_k})\To \h{A}\To \h{A'}\To 1.
\end{equation}
The connecting homomorphism $\delta$ sends any $U'\in (A')^{G_k}$ to
the homomorphism
$$\delta(U')\colon G_k\To \io,\quad \sigma\longmapsto
U\, ^{\sigma}\!U^{-1}, $$
where $U\in A$ is any preimage
of $U'$. From (\ref{b}) we get:

\begin{prop}\label{card}
The isotropy group $H_C$ of $C$ under the action of $H$ on $\tw$
corresponds to $\delta((A')^{G_k})$ under the isomorphism $H\simeq
\h{\io}$ of (\ref{hac}). Moreover, $$\tw=\coprod_{\xi'\in
\h{A'}}\op{Tw}^0(C_{\xi}/k), $$ where $\xi\in\h{A}$ is any choice of
a preimage of $\xi'$ and $C_{\xi}$ is the twist of $C$ corresponding
to $\xi$ by (\ref{tw}).
\end{prop}

\begin{pf} Since $\iota$ is central in $A$, the group $\h{\io}$ acts on
the set $\h{A}$ by multiplication and the orbits are the subsets of classes
having the same image in $\h{A'}$. On the other hand, it is easy to
check that the twists corresponding to the orbit of one class
$\xi\in\h{A}$ are precisely those in $\op{Tw}^0(C_{\xi}/k)$. \qed
\end{pf}

By the Proposition, the orbits of $\tw$ under the action of $H$ are
parameterized by $\h{A'}$ and the isotropy group of each $C_{\xi}$ is
$H_{C_{\xi}}\simeq\delta({A'_{\xi}}^{G_k})$, where $A'_{\xi}$ is the
reduced group of automorphisms of the curve $C_{\xi}$. These curves
$C_{\xi}$ may have in principle different isotropy subgroups. If $A$
is abelian, we get a more uniform description of $\tw$, since we have
then an isomorphism of pointed sets:
$$ \h{A}\simeq
\h{\io}/\delta((A')^{G_k})\times \h{A'},$$
and all curves $C_{\xi}$
have the same isotropy subgroup.

It is easy to describe $H_C$ in terms of a rational model of $C$ by
applying the results of section 1:

\begin{prop}\label{cqt}
Let $C$ be a curve of genus two defined over $k$ given by a rational
model with parameters $(a,b,c,d)$ as in section 1. Then,
\begin{enumerate}
\item If $C$ is supersingular, then $H_C=\delta_{abc}(\g_{abc})$ (see
Lemma \ref{delta}).
\item If $C$ is of type (1,1,1)-split, two at least of $a,b,c$ are equal
and the third coefficient (equal to the former two or not) does not
belong to $\as(k)$, then $H_C\simeq C_2$ is the subgroup generated by
this third coefficient.

\item If $C$ is of type (1,1,1)-quadratic, $b=0$ and $a\not \in\as(k)$,
then $H_C\simeq C_2$ is the subgroup generated by $a$.

\item In all other cases, $H_C$ is trivial.
\end{enumerate}
\end{prop}

\subsection{Twists of curves of genus two}
By Corollary \ref{jc},
each triple $j=(j_1,j_2,j_3)\in k^3$ can be identified with a unique
element in 
$\cc$, the quotient set of all
curves of genus two defined over $k$ classified up to isomorphism.
By Theorem \ref{igusa}, this set $\cc$ decomposes,
$$ \cc=\cc_{C_2}\cup
\cc_{C_2\times C_2}\cup \cc_{C_2\times S_3}\cup \cc_{M_{32}}\cup
\cc_{M_{160}}, $$
as the disjoint union of subfamilies gathering all
curves with isomorphic automorphism group. In this paragraph we find
an explicit description of these subfamilies and of the twists of
each curve in the family.

\subsubsection*{Curves with $\au\simeq C_2$}
This is the generic case. The set of values of the geometric
invariants is:
$$ j(\cc_{C_2})=J:=\{(j_1,j_2,j_3)\in k^3\tq j_1\ne
j_2j_3\}\sqcup\{(0,j_2,0)\in k^3\tq j_2\ne0\}. $$
We build up
$\cc_{C_2}$ by choosing for each $j\in J$ the curve with $j(C)=j$
indicated in the proof of Theorem \ref{cosdef}. For each curve
$C\in\cc_{C_2}$, $\tw=\twc$ is in bijection with the set of arithmetic
invariants $\h{A}\simeq \kas$.

If $k$ is the finite field with $q$ elements, there are $q^3-q^2+q-1$
values of the geometric invariants and the total number of curves
with $\au\simeq C_2$ is $2(q^3-q^2+q-1)$.

\subsubsection*{Curves with $\au\simeq C_2\times C_2$}
The set of values of the geometric invariants is:
$$
j(\cc_{C_2\times
C_2 })=J:=\{(j_1,j_2,j_3)\in k^3\tq j_1=j_2j_3,\ j_1\ne 0,\,j_2\ne
j_3^2\}.
$$
The condition $(j_1,j_2,j_3)\in J$ is equivalent to:
$$
x^3+j_3x^2+j_2x+j_1=(x+a)(x+c)^2,\quad a,c\in k^*,\, a\ne c. $$
Hence,
we can make the following choice for $\cc_{C_2\times C_2}$:
$$
\cc_{C_2\times C_2}=\{y^2+y=ax+\frac cx+\frac c{x+1}\,\big{|}\, a,c\in
k^*,\, a\ne c\}. $$
Let $C\in \cc_{C_2\times C_2}$ be fixed and let
$w\in\kb$ be such that $w^2+w=a$. The automorphism group of $C$ is
$A=\{1,\iota,U,\iota U\}$, where $U$ is the non-hyperelliptic
involution, $U(x,y)=(x+1,y+w)$. If $a\in\as(k)$, then $A$ has trivial
$G_k$-action. If $a\not\in \as(k)$, the automorphisms $U$ and $\iota
U$ are exchanged by the action of any $\sigma\in G_k$ with
non-trivial image in $\op{Gal}(k(w)/k)$. In both cases $A'$ has
trivial action and $H_C\simeq\delta(A')$ is the subgroup of $\kas$
generated by the class of $a$. Since $A$ is abelian,
\begin{myeqnarray*}
\lefteqn{\h{A}\simeq}\qquad\\&&
\simeq\h{\io}/\delta(A')\times \h{A'}\simeq \frac k{\{0,a\}+\as(k)}\times
\frac k{\as(k)}.
\end{myeqnarray*}%
A class $\xi'\in\h{A'}=\op{Hom}(G_k,A')$ is
identified with the element $r\in\kas$ representing the quadratic or
trivial extension of $k$ through which the homomorphism $\xi'$
factorizes. As a preimage $\xi\in\h{A}$ we can choose the class
corresponding to the twist:
$$ C_{\xi}: \qquad y^2+y=ax+\frac
c{x^2+x+r}\,, $$
which is isomorphic to $C$ via
\begin{equation}\label{isom}
(x,y)\mapsto(x+\theta,y+v),\quad \theta^2+\theta=r,\
v^2+v=a\theta.\end{equation}
We get the following parameterization of
$\tw$ in terms of the arithmetic invariants:
\begin{myeqnarray*}
\lefteqn{\tw=}\\&&
=\{y^2+y=ax+\frac
c{x^2+x+r}+d\,\big{|}\, (d,r)\in \frac k{\{0,a\}+\as(k)}\times \frac
k{\as(k)}\}.
\end{myeqnarray*}%

If $k$ is the finite field with $q$ elements, there are $(q-1)(q-2)$
values of the geometric invariants and the total number of curves
with $\au\simeq C_2\times C_2$ is:
$$ 4(\frac q2-1)(q-2)+2\frac q2(q-2)=(3q-4)(q-2).$$

\subsubsection*{Curves with $\au\simeq C_2\times S_3$}
The values of the geometric invariants and a concrete choice for
$\cc_{C_2\times S_3}$ are:
$$ j(\cc_{C_2\times S_3
})=\{(j_3^3,j_3^2,j_3)\in k^3\tq j_3\ne 0\},
$$
$$ \cc_{C_2\times
S_3}=\{y^2+y=ax+\frac ax+\frac a{x+1})\,\big{|}\, a\in k^*\}.
$$
Let
$C\in \cc_{C_2\times S_3}$ be fixed and let $w\in\kb$ be such that
$w^2+w=a$. The automorphisms
$$ U(x,y)=(x+1,y+w),\quad V(x,y)=(\frac 1{x+1},y), $$
generate a subgroup $S\subset A$, isomorphic to $S_3$,
such that $A=\io\times S$ (see the remark after Proposition
\ref{uuu}). Hence, $A'\simeq S_3$ and it is generated by the images $U',V'$
of $U,V$. If $a\in\as(k)$, then $A$ has
trivial $G_k$-action. If $a\not\in \as(k)$, only the automorphisms
$U$, $\iota U$ are not $G_k$-invariant; they are exchanged by the
action of all $\sigma\in G_k$ with non-trivial image in
$\op{Gal}(k(w)/k)$. In both cases $A'$ has trivial $G_k$ action, so
that:
$$ \h{A'}=\op{Hom}(G_k,A')\backslash \op{Inn}(A')\simeq \frac
k{\as(k)}\sqcup \ck.$$
Any non-trivial $r\in\kas$, associated to a
quadratic extension $K/k$, is identified with the class of the
1-cocycle determined by the isomorphism
$\op{Gal}(K,k)\stackrel{\sim}{\to}\{1,U'\}$. Any $s\in\ck$, associated
to a cubic extension $K/k$ with normal closure $\tilde{K}$, is
identified with the class of the 1-cocycle determined by an
isomorphism
$\op{Gal}(\tilde{K},k)\stackrel{\sim}{\to}\{1,V',(V')^2\}$, or
$\op{Gal}(\tilde{K},k)\stackrel{\sim}{\to}A'$. We can respectively
choose as twists $C_{\xi}$ associated to preimages $\xi\in\h{A}$:
$$
C_{\xi}\colon\qquad y^2+y=ax+\frac a{x^2+x+r},\quad y^2+y=\frac
{asx^2+asx+as(s+1)}{x^3+sx+s}.
$$
The curves $C_{\xi}$ with quadratic
denominator are isomorphic to $C$ via (\ref{isom}), whereas for those
with cubic denominator we can take:
$$ (x,y)\mapsto
(\frac{\theta(x+\theta)}{\theta'(x+\theta')},y+v),\quad
\theta^3+s\theta+s=0,\ v^2+v=aw. $$ By Propositions \ref{card} and
\ref{cqt}, $\tw$ is parameterized as the disjoint union:
\begin{myeqnarray*}
\lefteqn{\tw=}\\
&=&\{ y^2+y=ax+\frac
a{x^2+x+r}+d\,\big{|}\,(d,r)\in \frac k{\{0,a\}+\as(k)}\times \frac
k{\as(k)}\}\,\\
&&\sqcup \,\{y^2+y=\frac
{asx^2+asx+as(s+1)}{x^3+sx+s}+d,\big{|}\,(d,s)\in \frac
k{\as(k)}\times \ck\}.
\end{myeqnarray*}%
As mentioned in section 1, for certain
values of $s\in\ck$ one can use the curves, $y^2+y=as/(x^3+s)+d$,
instead of the above ones.

If $k$ is the finite field with $q$ elements, there are $q-1$ values
of the geometric invariants and the total number of curves with
$\au\simeq C_2\times S_3$ is:
$$2\frac q2+4(\frac q2-1)+2(q-1)=5q-6.$$

\subsubsection*{Curves with $\au\simeq M_{32}$}
The values of the geometric invariants and a concrete choice for
$\cc_{C_2\times M_{32}}$ are:
$$ j(\cc_{M_{32} })=\{(0,0,j_3)\in
k^3\tq j_3\ne 0\},$$
$$\cc_{M_{32}}=\{y^2+y=ax^5+ax^3\tq a\in k^*\}.
$$
Let $C\in \cc_{M_{32}}$ be fixed and denote by
$E(x):=E_{aa}(x)=a^4x^{16}+a^4x^8+a^2x^2+ax$ the linear separable
polynomial considered in Lemma \ref{delta}. Clearly, $\ke$ is a
subgroup of $\kb$ isomorphic to $C_2\times C_2\times C_2\times C_2$.
Recall that we have an exact sequence (\ref{meta}):
\begin{equation}\label{sex}
  1\To\io\To\au\To \ke\To 1.
\end{equation}
Given $\nu\in\ke$, the two automorphisms of $C$ lifting $\nu$ are
$$ \pm U_{\nu}(x,y)=(x+\nu,y+t_0+t_1x+t_2x^2),$$
where
$$
t_1=a\nu^2(1+\nu^2),\quad t_2=\sqrt{a\nu},\quad
t_0^2+t_0=a\nu^3(1+\nu^2).
$$
The plus or minus sign before $U_{\nu}$
correspond to any choice of the two possible values of $t_0$. Thus,
(\ref{sex}) is an exact sequence of $G_k$-groups and $A'\simeq \ke$
as $G_k$-groups.

From the exact sequence of abelian $G_k$-groups:
$$
1\To\ke\To\kb\stackrel{E}\To \kb\To 1, $$
we see that $\h{\ke}\simeq
k/E(k)$. For any $b\in k/E(k)$, the 1-cocycle associated to $b$ is
given by
$\{\xi'_{\sigma}=\beta+\,^{\sigma}\!\beta\}_{\sigma\in G_k}$,
where $\beta\in\kb$ satisfies $E(\beta)=b$.
Thus, we can choose the following twists associated to preimages in
$\h{A}$ of these cocycles:
$$C_{\xi}\colon\qquad
y^2+y=ax^5+bx^4+ax^3,\qquad b\in k/E(k).$$
Isomorphisms to $C$ are
given by $(x,y)\mapsto (x+\beta,y+t_0+t_1x+t_2x^2)$, where
$$
t_1=a\beta^2(1+\beta^2),\quad t_2=\sqrt{a\beta+b},\quad
t_0^2+t_0=a\beta^3(1+\beta^2)+b\beta^4. $$
By Propositions \ref{card}
and \ref{cqt},
\begin{myeqnarray*}
\lefteqn{\tw=}\quad\\&&=\coprod_{b\in k/E(k)}\{
y^2+y=ax^5+bx^4+ax^3+d\,\big{|}\,d\in \frac
k{\delta_{aba}(\g_{aba})+\as(k)}\}. 
\end{myeqnarray*}

If $k$ is the finite field with $q$ elements, there are $q-1$ values
of the geometric invariants and the total number of curves with
$\au\simeq M_{32}$ is $4q-5-[2]_{2|m}$.

\subsubsection*{Curves with $\au\simeq M_{160}$}
In this case $ \cc_{M_{160}}$ consists in one single curve $C$, with
model $y^2+y=x^5$ and $j$-invariant $j(C)=(0,0,0)$.

For any $a\in k^*$, we denote by $ E_a(x):=E_{a0}(x)=a^4x^{16}+ax$,
the $\fs$-linear separable polynomial considered in Lemma
\ref{delta}. Since $E_a(\la x)=\la E_a(x)$, for all $\la\in\mu_5$,
there is a well defined action of the group $\mu_5(k)$ on $k/E_a(k)
\simeq \h{\op{Ker}(E_a)}$.

For $a=1$, we have $\op{Ker}(E_1)=\fs$, and the exact sequence
(\ref{meta}) is:
\begin{equation}\label{set}
  1\To\io\To\au\To \mu_5\rtimes\fs\To 1.
\end{equation}
The two automorphisms of $C$ lifting any $(\la,\nu)\in
\mu_5\rtimes\fs$ are
$$ \pm U_{(\la,\nu)}(x,y)=(\la
x+\nu,y+t_0+t_1x+t_2x^2),$$
where
$$ t_1=\la\nu^4,\quad
t_2=\la^2\nu^8,\quad t_0^2+t_0=\nu^5. $$
The plus or minus sign
before $U_{(\la,\nu)}$ corresponds to any choice of the two possible
values of $t_0$. Thus, (\ref{set}) is an exact sequence of
$G_k$-groups and $A'\simeq \mu_5\rtimes \fs$ as $G_k$-groups.

From the split exact sequence of $G_k$-groups:
$$ 1\To\fs\To A'\To \mu_5\To 1, $$
we get an exact sequence of pointed sets:
\begin{equation}\label{ult}
\h{\fs}\stackrel{i}\To\h{A'}\stackrel{\pi}\To \h{\mu_5}\To 1.
\end{equation}
Under the identification $\h{\fs}\simeq k/E_1(k)$, we have that 
$i(b)=i(b')$ if and only if there exists $\la\in\mu_5(k)$ such that
$b'=\la b$; hence, we get from (\ref{ult}) an exact sequence of
pointed sets:
$$ 1\To(k/E_1(k))\backslash
\mu_5(k)\stackrel{i}\To\h{A'}\stackrel{\pi}\To \h{\mu_5}\To 1.$$
In a
similar way, for any other value of $a\in k^*$, we can consider the mapping
$$
(k/E_a(k))\backslash \mu_5(k)\To \h{A'},
$$
sending the class of any
$b\in k$ to the class of the 1-cocycle:
$$
\{\xi'_{\sigma}=(\alpha\,^{\sigma}\!
\alpha^{-1},\alpha(\beta+\,^{\sigma}\!\beta))
\}_{\sigma\in G_k}, 
$$
where $\alpha,\beta\in
\kb$ satisfy $\alpha^5=a$ and $E_a(\beta)=b$. It is not difficult to
check that this mapping is 1-1 with image $\pi^{-1}([a])$, where
$[a]\in k^*/(k^*)^5\simeq \h{\mu_5}$ is the class of $a$. Thus,
$$H^1(G_k,A')\simeq \coprod_{a\in k^*/(k^*)^5}
(k/E_a(k))\backslash \mu_5(k).$$
We can choose as twists associated to
preimages in $\h{A}$:
$$C_{\xi}\colon \qquad y^2+y=a x^5+b x^4,\quad
a\in k^*/(k^*)^5,\ b\in (k/E_a(k))\backslash \mu_5(k),$$
isomorphic to
$C$ via $(x,y)\mapsto (\la(x+\beta),y+t_0+t_1x+t_2x^2)$, where
$$
\la^5=a,\quad t_1=a\beta^4,\quad t_2=a^2\beta^8,\quad
t_0^2+t_0=a\beta^5.
$$
By Propositions \ref{card} and \ref{cqt}, we can
parameterize $\tw$ as:
\begin{myeqnarray*}
  \lefteqn{\tw=
  \coprod_{a\in k^*/(k^*)^5}\,\, 
\coprod_{b\in (k/E_a(k))\backslash \mu_5(k)}}\qquad\qquad\qquad\qquad\\
&&\{y^2+y=ax^5+bx^4+d\,\big{|}\,d\in \frac
k{\delta_{ab0}(\g_{ab0})+\as(k)}\}.
\end{myeqnarray*}

If $k$ is the finite field with $q=2^m$ elements, there is only one
value of the geometric invariants and the total number of curves with
$\au\simeq M_{160}$ is $3+[2]_{2|m}+[8]_{4|m}$.

In fact, if $4\nmid m
$, the groups $k^*/(k^*)^5$, $\mu_5(k)$  vanish and
$$
H^1(G_k,A')\simeq k/E_1(k)=k/\as^4(k) $$
has 2 or 4 elements
according to $m$ odd or even. The twists of $C$ are:
$$ \tw=\{y^2+y=x^5+\epsilon
x^4\tq \epsilon\in\mathbb{F}\}\cup \{y^2+y=x^5+x^4+d_0\}, $$
where
$\mathbb{F}=\f2,d_0=1$ if $m$ is odd, whereas for $m$ even,
$\mathbb{F}=\mathbb{F}_4$ and $d_0$ is any choice of one of the two
elements in $\mathbb{F}_4-\f2$.

If $4|m$, we have $k/E_a(k)=0$ for all $a\not\in(k^*)^5$ and
$(k/\as^4(k))\backslash \mu_5$ has 4 elements. Hence, $H^1(G_k,A')$
has 8 elements and the 13 twists of $C$ are:
\begin{myeqnarray*}
\tw
&=&\{y^2+y=ax^5+d,\tq (a,d)\in k^*/(k^*)^5\times \kas\}\,\cup\\
&&{}\cup\{y^2+y=x^5+bx^4,\tq b\in ((k/\as^4(k))-\{0\})\backslash
\mu_5(k)\}.
\end{myeqnarray*}

\end{article}

\end{document}